\documentclass[final,1p,number,sort&compress]{siamart1116}

\usepackage[english]{babel}
\usepackage{graphicx,epstopdf,epsfig}
\usepackage{amsfonts,epsfig,fancyhdr,graphics, hyperref,amsmath,amssymb}

\numberwithin{equation}{section}
\usepackage{mathtools,mathdots}
\usepackage{algorithm}
\usepackage{algpseudocode}
\usepackage{amssymb}
\usepackage{amsmath}
\usepackage{booktabs,url}

\newcommand{\Names}{Heike Fa\ss bender, Martin Halwa\ss }
\newcommand{\Title}{A note on the SVD of idempotent and involutory matrices}

\newtheorem{remark}[theorem]{Remark}

\renewcommand{\theequation}{\thesection.\arabic{equation}}
\renewtheorem{theorem}{Theorem}[section]

\numberwithin{equation}{section}
\begin{document}

\setcounter{page}{1}

\thispagestyle{empty}

%Insert the title of the paper
 \title{A note on the singular value decomposition of idempotent  and involutory matrices}

\author{Heike Fa\ss bender\thanks{Institute for Numerical Analysis, TU Braunschweig,
Braunschweig, Germany (h.fassbender@tu-braunschweig.de).}
% Remember to put \and between any two authors
\and
{
Martin Halwa\ss \thanks{Loitz, Germany (martin.halwass@web.de).}
}}

\pagestyle{myheadings}
\markboth{\Names}{\Title}

\maketitle

\begin{abstract}
It is known  that singular values of idempotent matrices are either zero or larger or equal to one \cite{HouC63}. We state exactly how many
 singular values greater than one, equal to one, and equal to zero  there are.
Moreover, we derive a singular value decomposition of idempotent matrices which reveals a tight relationship between its left and right singular vectors. The same idea is used to augment a discovery regarding the singular values of involutory matrices as presented in \cite{FasH20}.
\end{abstract}

\begin{keywords}
idempotent matrix, singular value decomposition, QR decomposition, Schur decomposition, involutory matrix
\end{keywords}

\begin{AMS}
15A24 $\cdot$ 65F15
\end{AMS}

%%%%%%%%%%%%%%%%%%%%%%%%%%%%%%%%%%%%%%%

%%%%%%%%%%%%%%%%%%%%%%%%%%%%%%%%%%%%%%%%%%%%%%%%%%%%%%%%%%%%%
\section{Introduction}
An idempotent matrix $M \in \mathbb{C}^{n \times n}$  is a matrix which, when multiplied by itself, yields itself. That is, $M^2 = M.$ It is well known that an idempotent matrix is always diagonalizable \cite[Problem 3.3.P3]{HorJ13}. Its eigenvalues are either $0$ or $1.$
Thus, idempotent matrices are of the form $M = ZDZ^{-1},$
where D is a diagonal matrix having only 0’s and 1’s on the diagonal, and Z is regular.
 Hence, its determinant is either  $1$ in case $M=D=I$ and $0$ otherwise.
Householder and Carpenter noted in  \cite{HouC63} that singular values of $M$ are either zero or larger or equal to one. Clearly, the number of singular values equal to zero corresponds to $n-\operatorname{rank}(M).$ In \cite{HouC63}, however, nothing is said about the number of singular values larger than one and the number of those equal to one.

In this note, we will shed some more light on the singular value decomposition (SVD) of idempotent matrices. In particular, we derive an SVD of $M$ from which the number of singular values larger than one, equal to one and equal to zero can be read off in concrete terms.
 In particular, making use of a result from  \cite{Dok91}, we prove that there are exactly $t = \operatorname{rank}(M) -\operatorname{dim}(\operatorname{null}(I-MM^H))$ singular values larger than one and $\operatorname{rank}(M) -t= \operatorname{dim}(\operatorname{null}(I-MM^H))$ singular values equal to one. Moreover, we will see a close connection between the left and right singular vectors.
Finally, in Section \ref{sec4}, the result from \cite{Dok91}  is used in order to complete a finding on the singular values of involutory matrices given in \cite{FasH20}.
%%%%%%%%%%%%%%%%%%%%%%%%%%%%%%%%%%%%%%%%%%%%%%%%%%%%
%%%%%%%%%%%%%%%%%%%%%%%%%%%%%%%%%%%%%%%%%%%%%%%%%%%%

\section{Singular value decomposition of idempotent matrices}\label{sec2}
 First, we repeat some basic properties of an SVD (not only of an idempotent matrix). Then, we consider the SVD of an idempotent matrix and show that the singular values $\sigma_j$ are given by the inverse of the cosine of the angle $\psi_j$ between the subspaces spanned by the corresponding singular left and right vectors. Next, we review some results on idempotent matrices from the literature which give us a condensed Schur form $N$ of $M$ in which the only nonzero entries that are not on the diagonal are $\tan \psi_j.$ From the SVD of $N$, we obtain the SVD of $M$ which reveals the number of singular values larger than one, equal to one and equal to zero as well as a close connection between the left and right singular vectors.

Let $M \in \mathbb{C}^{n \times n}$ be an idempotent matrix of rank $r.$ We are interested in its SVD $U\Sigma V^H,$ where $U \in \mathbb{C}^{n\times n}$ and $V\in \mathbb{C}^{n\times n}$ are unitary and $\Sigma = \left[\begin{smallmatrix}\Omega&\\&0\end{smallmatrix}\right] \in \mathbb{R}^{n\times n}$ where $\Omega \in \mathbb{R}^{r \times r}$ is a diagonal matrix whose nonnegative diagonal elements $\sigma_j$ are ordered in decreasing order, $\sigma_1 \geq \sigma_2 \geq \cdots \geq \sigma_r >0.$
It is well-known that the columns of $U$ are eigenvectors of $MM^H$ as $MM^H= U\Sigma^2U^H$ and that the columns of $V$ are eigenvectors of $M^HM$ as $M^HM= V\Sigma^2V^H.$ Moreover, the nonzero elements of $\Sigma$ (the nonzero singular values) are the square roots of the nonzero eigenvalues of $MM^H$ as well as $M^HM$.
The left singular vectors (= column vectors of $U$) corresponding to the nonzero singular values of $M$ span the range of $M.$ Hence, the right singular vectors (= column vectors of $V$) corresponding to the nonzero singular values of $M$ span the range of $M^H.$

From $M^2 = U\Sigma V^HU\Sigma V^H = U \Sigma V^H = M$ we have $\Sigma V^HU\Sigma =  \Sigma.$ That is, with $U =\left[U_r ~~U_s\right], U_r \in \mathbb{C}^{n \times r}$ and $V =\left[V_r ~~V_s\right], V_r \in \mathbb{C}^{n \times r}$ we have
\[
\Sigma V^HU\Sigma =  \begin{bmatrix}\Omega V_r^H U_r \Omega&\\&0_{s\times s}\end{bmatrix} = \begin{bmatrix}\Omega&\\&0_{s\times s}\end{bmatrix}=\Sigma
\]
 as  $U\Sigma =\left[ U_r\Omega ~~0_{n\times s}\right]$ and $V\Sigma=\left[V_r\Omega ~~0_{n\times s}\right].$ Therefore, $ V_r^H U_r^H = \Omega^{-1}.$
 Hence, the following lemma holds.
\begin{lemma}
Let $M\in \mathbb{C}^{n \times n}$ be an idempotent matrix of rank $r$ with the SVD $M = U \Sigma V^H.$
Let $u_j$ and $v_j$ be the left and right singular vectors, respectively, of $M$ associated with the nonzero singular values $\sigma_j, j = 1, 2, \ldots, r.$ Then, $u_j^Hv_j = \sigma_j^{-1}.$
\end{lemma}
As $u_j$ and $v_j$ are vectors of length $1,$ we obtain from the definition of the angle between the subspace $\{u_j\}$ and $\{v_j\}$  \cite[Problem 5.1.P3]{HorJ13} (see \cite[Formuale (23)]{BehS94} for a similar derivation) that
\[
|u_j^Hv_j |= \cos \psi_j = \sigma_j^{-1} \neq 0
\]
holds with
\begin{equation}\label{eq_psi}
\frac{\pi}{2} > \psi_1 \geq \psi_2 \geq \cdots \geq \psi_t > \psi_{t+1}=\cdots = \psi_r = 0.
\end{equation}

In \cite[Theorem 1]{Dok91}, it is shown that every idempotent matrix is unitarily similar to a real $n \times n$ block diagonal matrix
\begin{equation}\label{bd}
\widehat{N}= \begin{bmatrix} 1 &\tau_1\\0&0\end{bmatrix} \oplus \cdots \oplus \begin{bmatrix} 1 &\tau_t\\0&0\end{bmatrix} \oplus I_{r-t} \oplus 0_{s-t}
\end{equation}
where $\tau_1 \geq \tau_2 \geq \cdots \geq \tau_t,$ $r = \operatorname{rank}(M),$ $s = n-r,$ and
\[
t = \operatorname{rank}(M)-\dim(\operatorname{null}(I-MM^H)).
\]
In \cite{Ikr96} it is shown that the $\tau_j$'s  in \eqref{bd} are the tangents of the principal angles between $\operatorname{range}(M)$ and $\operatorname{range}(M^H).$
Thus, in \cite{Ikr96} it is observed  that $\tau _j = \tan \psi_j$ for $j = 1, \ldots, t$ holds for the angles $\psi_j$ as in \eqref{eq_psi}. The connection between the $\tau_j$ in \eqref{bd} and the singular values of $M$ has also been observed in \cite[Corollary 3.4.3.3.]{HorJ13}. There it is noted that $\tau_j = \sqrt{\sigma_j^2-1}$  holds. Inserting $\tau_j = \tan \psi_j,$ this gives $\sigma_j = 1/\cos \psi_j.$ Moreover,  in \cite[Theorem 58]{Gal08}  it is noted that an idempotent matrix $M$  is unitarily similar to a matrix of the form \eqref{bd} where $\tau_j = \cot \theta_j$ for certain principal angles between two subspaces associated with $M.$

For our purposes, we will consider a permuted version of the matrix $\widehat{N}$,
\begin{equation}\label{eq:N}
N = \begin{bmatrix}I_t &   & & T\\
& I_{r-t}& \\
 && 0_{s-t}&\\
&&&0_{t}
\end{bmatrix}
\end{equation}
with
$T = \operatorname{diag}(\tan \psi_1, \tan \psi_2, \ldots, \tan \psi_t).$
Thus, from \eqref{bd} we have a Schur decomposition of $M$
\begin{equation}\label{schur}
M = \mathcal{U} N \mathcal{U}^H
\end{equation}
with a unitary matrix $\mathcal{U} \in \mathbb{C}^{n \times n}.$

Next, the SVD of $N$ is described. With this, we will be able to give the SVD of $M$ revealing the desired information.
\begin{lemma}[SVD of $N$]\label{svdN}
The real $n \times n$ matrix $N$  \eqref{eq:N} has the SVD
\[
N = \begin{bmatrix}I_t & &&T \\
& I_{r-t}& \\
 && 0_{s-t}&\\
&&&0_{t}
\end{bmatrix} = \mathcal{S}O^T
\]
with the $n \times n$ real diagonal matrix
\[
\mathcal{S} = \begin{bmatrix} C^{-1}& \\ & I_{r-t}&\\ && 0_{s-t}\\&&&0_{t} \end{bmatrix}
\]
and an $n \times n$ real orthogonal matrix
\[
O =  \begin{bmatrix} C& &&ES\\  && I_{n-2t}\\S&&& -EC  \end{bmatrix}
\]
with $C = \operatorname{diag}(\cos \psi_1, \cos \psi_2, \ldots, \cos \psi_t),$ $S = \operatorname{diag}(\sin \psi_1, \sin \psi_2, \ldots, \sin \psi_t),$ and $ E= \operatorname{diag}(\pm 1, \ldots, \pm 1) \in \mathbb{R}^{t \times t}.$
\end{lemma}
\begin{proof}
We determine the SVD of $N$ from the eigendecompositions of $NN^T$ and $N^TN.$ First, form $NN^T$ and use $\cos^{-2}\psi = 1 + \tan^2 \psi$ to obtain
\begin{align*}
NN^T &= \begin{bmatrix} I_t+T^2 & \\& I_{r-t}&\\&&0_{s-t}\\&&&0_{t} \end{bmatrix}
= \begin{bmatrix} C^{-2} & \\& I_{r-t}&\\&&0_{s-t}\\&&&0_{t} \end{bmatrix}=\mathcal{S}^2.
\end{align*}
Next, we form the symmetric matrix $N^TN$
\[
N^TN = \begin{bmatrix} I_t & &&T\\ &I_{r-t}\\  &&0_{s-t}\\ T & &&T^2 \\ \end{bmatrix}.
\]
Each of the four blocks in $\left[ \begin{smallmatrix} I_t & T\\ T & T^2  \end{smallmatrix}\right]$ is diagonal.  In order to compute an eigendecomposition of $N^TN$ we need to eliminate the diagonal entries in the (2,1)-block. This can be done without introducing any new entries in the (2,1)-, the (1,1)- and the (2,2)-block) by observing that
\[
\begin{bmatrix} 1 & \tan \psi_i \\ \tan \psi_i & \tan^2 \psi_i\end{bmatrix}=
\begin{bmatrix} \cos \psi_i & -\sin \psi_i \\ \sin \psi_i & \cos \psi_i\end{bmatrix}
\begin{bmatrix} 1 +\tan^2 \psi_i& 0 \\ 0 & 0\end{bmatrix}\begin{bmatrix} \cos \psi_i & \sin \psi_i \\ -\sin \psi_i & \cos \psi_i\end{bmatrix}
\]
as well as
\[
\begin{bmatrix} 1 & \tan \psi_i \\ \tan \psi_i & \tan^2 \psi_i\end{bmatrix}=
\begin{bmatrix} \cos \psi_i & \sin \psi_i \\ \sin \psi_i & -\cos \psi_i\end{bmatrix}
\begin{bmatrix} 1 +\tan^2 \psi_i& 0 \\ 0 & 0\end{bmatrix}\begin{bmatrix} \cos \psi_i & \sin \psi_i \\ \sin \psi_i & -\cos \psi_i\end{bmatrix}
\]
holds. Thus, we have $N^TN = O\mathcal{S}^2O^T.$ Hence, we see that matrix of the left singular vectors of $N$ is given by the identity and the matrix of the right singular vector is given by $O,$ while the singular values are given by the diagonal elements of $\mathcal{S}.$ This proves the statement of the lemma.
\end{proof}

\begin{theorem}[SVD of idempotent matrices]
Let $M \in \mathbb{C}^{n \times n}$ be an idempotent matrix of rank $r.$ Then, the SVD of $M$ is given by
\begin{equation}\label{svdC2}
M =  \mathcal{U} \mathcal{S} O^T\mathcal{U}^H = \mathcal{U} \mathcal{S} \mathcal{V}^H
\end{equation}
where $\mathcal{U}  \in \mathbb{C}^{n \times n}$ and $\mathcal{V} = \mathcal{U}O \in \mathbb{C}^{n \times n}$ are unitary.  The diagonal matrix $\mathcal{S} \in  \mathbb{R}^{n \times n}$ and the orthogonal matrix $O\in \mathbb{R}^{n \times n}$ are as in Lemma \ref{svdN}.
\end{theorem}
The SVD \eqref{svdC2} reveals that for every idempotent matrix $M,$ besides the $n-r$ singular values equal to $0,$ there are $t$ singular values greater than $1$ and $r-t$ singular values equal to $1,$ where $t= \operatorname{rank}(M) -\operatorname{dim}(\operatorname{null}(I-MM^H)).$ Moreover, there is a close connection between the left and the right singular vectors, as $\mathcal{V} = \mathcal{U}O $ holds.

\begin{remark}
Making use of $\mathcal{U}=\mathcal{V}O^T$ in \eqref{svdC2} we obtain $M=\mathcal{V}O^T\mathcal{S}\mathcal{V}^H.$ As
\[
O^T\mathcal{S} = \begin{bmatrix} I_{n-t}&\\&E\end{bmatrix} N^T  \begin{bmatrix} I_{n-t}&\\&E\end{bmatrix},
\]
this leads to the SVD of $M$ in the form $M = \mathcal{W}O\mathcal{S}\mathcal{W}^H$ with the unitary matrix $\mathcal{W}=\mathcal{V}(I_{n-t}\oplus E).$
\end{remark}

\begin{remark}\label{rem1}
We can say a little bit more about the unitary matrix $\mathcal{U}$ by making use of the fact that  the columns of $\mathcal{U}$ are eigenvectors of $MM^H.$
As a first step, a suitable rank-revealing QR decomposition \cite{GolvL13} of $M$ is determined. For this, a permutation matrix $\Pi$  is chosen so that the first $r$ columns of $A=\Pi^TM\Pi$ are linearly independent.  The resulting matrix  $A$ is idempotent. The permutation $\Pi$ is in addition chosen such that the diagonal elements of the triangular factor of the rank-revealing QR decomposition form a monotonically  decreasing sequence of positive real numbers, see, e.g., \cite[Page 12, Theorem (3.15)]{LawH74} and \cite[Page 283]{KieS88}. That is,  we have
\begin{align}
A &=\Pi^TM\Pi = QR = \begin{bmatrix}Q_1 \mid Q_2\end{bmatrix} \begin{bmatrix} \begin{array}{c}R_1 \\ \hline  0_{s\times n} \end{array}\end{bmatrix}
= \begin{bmatrix}Q_1 \mid Q_2\end{bmatrix} \begin{bmatrix} \begin{array}{cccc|cccc} \mathfrak{r}_{11} &*& \cdots &* & * &\cdots & *\\ &\ddots&\ddots& \vdots & \vdots & & \vdots\\&& \mathfrak{r}_{r-1,r-1}&*& * &\cdots & *\\ &&&\mathfrak{r}_{rr} & *& \cdots & *\\ \hline
\multicolumn{4}{c|}{0_{s\times r} } & \multicolumn{3}{c}{0_{s\times s}}\end{array}  \end{bmatrix}, \label{qr_A}
\end{align}
with an unitary matrix $Q \in \mathbb{C}^{n \times n},$ orthonormal matrices $Q_1 \in \mathbb{C}^{n \times r} $ and $Q_2 \in \mathbb{C}^{n \times s}, $ an upper triangular matrix $R_1 \in \mathbb{C}^{r \times n}$ and $\mathfrak{r}_{11} \geq \cdots \geq \mathfrak{r}_{r-1,r-1} \geq \mathfrak{r}_{rr} > 0.$ $Q_1$ and $R_1$ are uniquely determined, whereas the columns of $Q_2$ are not uniquely determined. In particular, any order of those columns is possible.

Next, we derive an eigendecomposition of $AA^H$. In order to do so, first consider a unitary similarity transformation of $A$  with the  unitary factor of the QR decomposition \eqref{qr_A}
\begin{equation}
Q^HAQ = RQ  = \begin{bmatrix}R_1Q_1 & R_1Q_2\\0_{s\times r}&0_{s\times s}\end{bmatrix} = \begin{bmatrix}I_r & X\\0_{s\times r}&0_{s\times s}\end{bmatrix},\label{eq_simqrA}
\end{equation}
with $X = R_1Q_2 \in \mathbb{C}^{r \times s}$.  The identity $R_1Q_1 = I_r $ follows due to the idempotency of $A$: multiplying $A = Q_1R_1 = A^2 = Q_1R_1Q_1R_1$ from the left by $Q_1^H$ and from the right by $R_1^H$ yields $R_1R_1^H = R_1Q_1R_1R_1^H.$
As $R_1$ is off full row rank, $R_1R_1^H$ is a positive definite Hermitian matrix. Hence, its inverse exists. Thus, from the above equation, we have $R_1Q_1=I_r.$ This property was already noted (in a slightly more general fashion) as formula (3)  in \cite{HouC63}.

Next, construct $Q^HAA^HQ$
\begin{equation} \label{eq_AAH}
Q^HAA^HQ = \begin{bmatrix}I_r & X\\0_{s\times r}&0_{s\times s}\end{bmatrix}\begin{bmatrix}I_r & 0_{r\times s}\\X^H&0_{s\times s}\end{bmatrix} = \begin{bmatrix}XX^H+I_r & 0_{r\times s}\\0_{s\times r}&0_{s\times s}\end{bmatrix}.
\end{equation}
By construction, $XX^H+I_r$ is a positive definite Hermitian matrix. Denote its eigendecomposition by
\begin{equation}\label{eq_U1}
XX^H+I_r = U_1\Lambda U_1^H
\end{equation}
 where $\Lambda \in \mathbb{R}^{r\times r}$ is a diagonal matrix with positive diagonal entries and $U_1 \in \mathbb{C}^{r \times r}$ is unitary. Assume that $U_1$ is chosen such that the eigenvalues appear in decreasing order on the respective diagonal (note, that even so, $U_1$ does not need to be unique in case of multiple eigenvalues). Substituting this into \eqref{eq_AAH} and rearranging yields the eigendecomposition
\[
AA^H = Q\begin{bmatrix}U_1& \\ &I_s\end{bmatrix}\begin{bmatrix}\Lambda& \\ &0_{s\times s}\end{bmatrix}\begin{bmatrix}U_1 & \\ &I_s\end{bmatrix}^H Q^H.
\]
Since $AA^H = \Pi^T MM^H\Pi$, we can read off parts of the SVD of $M.$ First, we see that $\sqrt{\Lambda} =\left[\begin{smallmatrix} C^{-1}& \\ &I_{r-t}\end{smallmatrix}\right]$ must hold. Thus, the leading $r \times r$ nonzero diagonal block of matrix $\mathcal{S}$ in the SVD of $M$ equals $\sqrt{\Lambda},$ $\mathcal{S} = \left[\begin{smallmatrix} \sqrt{\Lambda}&\\& 0_{s \times s}\end{smallmatrix}\right].$
Moreover, $\mathcal{U}$ can be expressed as
\[
\mathcal{U} = \Pi Q \begin{bmatrix}U_1& \\ & U_2\end{bmatrix}
\]
with the unitary matrix $U_1  \in \mathbb{C}^{r \times r}$ as in \eqref{eq_U1}, an arbitrary unitary matrix $U_2 \in \mathbb{C}^{s \times s}$ and $Q, \Pi$ from \eqref{qr_A}.

By the way, we see from \eqref{eq_U1} that $\Lambda-I_r = \left[\begin{smallmatrix} C^{-2}-I_r&\\&0_{r-t}\end{smallmatrix}\right] = \left[\begin{smallmatrix} T^{2}&\\&0_{r-t}\end{smallmatrix}\right] =U_1^HXX^HU_1$ holds. Thus, the nonzero singular values of $X$ are given by the diagonal elements of $T.$
\end{remark}
%%%%%%%%%%%%%%%%%%%%%%%%%%%%%%%%%%%%%%%%%%%%%%%%%%%%%%%%%
%%%%%%%%%%%%%%%%%%%%%%%%%%%%%%%%%%%%%%%%%%%%%%%%%%%%%%%%%
%%%%%%%%%%%%%%%%%%%%%%%%%%%%%%%%%%%%%%%%%%%%%%%%%%%%%%%%%
\section{An SVD of involutory matrices}\label{sec4}
For any idempotent matrix $M$, the  matrix
\[
B = \pm (2M-I)
\]
is involutory, that is, $B^2 = I.$ The approach of the previous section can be used to derive an SVD of $B.$  To simplify notation, we will consider only the case $B = 2M-I.$
First, note that with \eqref{eq_psi} and $\nu =\min\{r,s\}$  we can write \eqref{bd} as
\begin{equation}\label{bd1}
\widehat{N}= \left\{ \begin{array}{cc}
\begin{bmatrix} 1 &\tan \psi_1\\0&0\end{bmatrix} \oplus \cdots \oplus \begin{bmatrix} 1 &\tan \psi_t\\0&0\end{bmatrix} \oplus
\begin{bmatrix} 1 &0\\0&0\end{bmatrix} \oplus \cdots \oplus \begin{bmatrix} 1 &0 \\0&0\end{bmatrix} \oplus 0_{s-r} & \qquad r \leq s\\
\begin{bmatrix} 1 &\tan \psi_1\\0&0\end{bmatrix} \oplus \cdots \oplus \begin{bmatrix} 1 &\tan \psi_t\\0&0\end{bmatrix} \oplus
\begin{bmatrix} 1 &0\\0&0\end{bmatrix} \oplus \cdots \oplus \begin{bmatrix} 1 &0 \\0&0\end{bmatrix} \oplus I_{r-s} & \qquad r \geq s
\end{array} \right.
\end{equation}
with $\nu-t$ blocks of the form $\left[\begin{smallmatrix}1 & \tan \psi_j\\\ 0 & 0\end{smallmatrix}\right] = \left[\begin{smallmatrix}1 & 0\\ 0 & 0\end{smallmatrix}\right], j = t+1, \ldots, \nu.$
Instead of \eqref{eq:N}, we will consider the permuted matrix
\begin{equation}\label{eq:N1}
N = \begin{bmatrix}I_\nu &  & T\\
& \Upsilon & \\
&&0_\nu
\end{bmatrix}
\end{equation}
with
$T = \operatorname{diag}(\tan \psi_1, \tan \psi_2, \ldots, \tan \psi_t, 0, \ldots, 0) \in \mathbb{R}^{\nu \times \nu}$ and $\Upsilon =0_{n-2\nu}$ in case $r\leq s$ and $\Upsilon =I_{n-2\nu}$ in case $r \geq s.$
Thus, from  \cite[Theorem 1]{Dok91} we have a Schur decomposition of $M$
\begin{equation}\label{schurx}
M = \mathcal{U} N \mathcal{U}^H
\end{equation}
with a unitary matrix $\mathcal{U} \in \mathbb{C}^{n \times n}.$

With \eqref{schurx} we have a condensed Schur form of $B,$
\[
B = (2 \mathcal{U}N\mathcal{U}^H -I)=  \mathcal{U}(2N-I)\mathcal{U}^H.
\]
Let
\[
\mathcal{N} =2N-I = \begin{bmatrix} I_\nu &&2T \\ &\pm I_{n-2\nu}& \\&&-I_{\nu} \end{bmatrix}
\]
where (here and in the rest of this section)  $-I_{n-2\nu}$ is used in case $r \leq s$ and $I_{n-2\nu}$ otherwise.
The SVD of $\mathcal{N}$ can be deduced from the observation that
\[
\begin{bmatrix} 1 & 2\tan \psi\\ & -1\end{bmatrix} =
\begin{bmatrix} \sin \phi & \cos \phi\\ -\cos \phi & \sin \phi\end{bmatrix}
\begin{bmatrix} \tan \phi & \\  & \tan^{-1} \phi\end{bmatrix}
\begin{bmatrix} \cos \phi & \sin \phi\\ \sin \phi & -\cos \phi\end{bmatrix}^T
\]
holds with $\phi = \frac{1}{2} \left(\frac{\pi}{2} + \psi\right).$
Thus, $\mathcal{N}$ has the SVD
\[
\mathcal{N} =
\begin{bmatrix} S_\mathcal{N}  && C_\mathcal{N}  \\ & \pm I_{n-2\nu}\\  -C_\mathcal{N}  && S_\mathcal{N}   \\ \end{bmatrix}
\begin{bmatrix} \Sigma_\mathcal{N}  & \\ & I_{n-2\nu}\\&&& \Sigma_\mathcal{N} ^{-1}\end{bmatrix}
\begin{bmatrix} C_\mathcal{N}  && S_\mathcal{N}    \\ &  I_{n-2\nu}\\ S_\mathcal{N}  && -C_\mathcal{N}  \end{bmatrix}
\]
with
\begin{align*}
C_\mathcal{N}  &= \operatorname{diag}(\cos \phi_1, \ldots, \cos \phi_t, \cos \phi_{t+1},\ldots, \cos \phi_\nu), \\
S_\mathcal{N}  &= \operatorname{diag}(\sin \phi_1, \ldots, \sin \phi_t, \sin \phi_{t+1},\ldots, \sin \phi_\nu), \\
 \Sigma_\mathcal{N}  &= \operatorname{diag}(\tan \phi_1, \ldots, \tan \phi_t, \tan \phi_{t+1},\ldots, \tan \phi_\nu),
\end{align*}
for  $\phi_j = \frac{1}{2} \left(\frac{\pi}{2} + \psi_j\right), j = 1, \ldots, t$ and  $\phi_j = \frac{\pi}{4},  j = t+1, \ldots, \nu.$

Denote the SVD of $\mathcal{N}$ by $\mathcal{N} = U_\mathcal{N}\mathfrak{S}_\mathcal{N}V_\mathcal{N}^T$ and note that
\[
U_\mathcal{N} = V_\mathcal{N} T_\mathcal{N} \quad \text{ and } \quad V_\mathcal{N} = U_\mathcal{N} T_\mathcal{N}
\]
holds for
\[
T_\mathcal{N}  = \begin{bmatrix} &&I_\nu & \\  &\pm I_{n-2\nu}\\I_\nu \end{bmatrix}.
\]
Thus, the left singular vectors of $\mathcal{N}$ result from its right singular vectors by permutation.
In summary, we have the SVD of $B =  \mathcal{U}\mathcal{N}\mathcal{U}^H$
\[
B =  \underbrace{\mathcal{U}U_\mathcal{N}}_{\widetilde U_{1}}\mathfrak{S}_\mathcal{N}\underbrace{T_\mathcal{N}U_\mathcal{N}^T \mathcal{U}^H}_{\widetilde V_{1}^H} = \underbrace{\mathcal{U}V_\mathcal{N}T_\mathcal{N}}_{\widetilde U_{2}}\mathfrak{S}_\mathcal{N}\underbrace{V_\mathcal{N}^T \mathcal{U}^H}_{\widetilde V_{2}^H}
\]
with unitary matrices
\[
\widetilde V_{1} = \widetilde U_{1}T_\mathcal{N} \quad \text{ and } \quad \widetilde U_{2} = \widetilde V_{2}T_\mathcal{N}.
\]
From $\mathfrak{S}_\mathcal{N},$ we see that $2\nu$ singular values appear in pairs $(\sigma, \frac{1}{\sigma}).$ Their left and right singular vectors are related. Among the singular values which appear in pairs $(\sigma, \frac{1}{\sigma})$  there are $\nu-t$  pairs $(1,1),$ while for the other $t$ pairs $(\sigma, \frac{1}{\sigma})\neq (1,1)$ holds.
Moreover, there are $n-2\nu$ singular values of size $1.$
The same statement has been derived in \cite[Theorem 3]{FasH20} by different means. What is new here, is that we can now specify $t$ in concrete terms, $t = \operatorname{rank}(M)-\dim(\operatorname{null}(I-MM^H))$ for the idempotent matrix $M$ which defines the involutory matrix $B.$ Moreover, we see that the singular values $\sigma$ not equal to $1$ are related to the singular values of the corresponding idempotent matrix as $ \sigma_j = \tan \phi_j$ holds with $\phi_j = \frac{1}{2}(\frac{\pi}{2}+\psi_j), j = 1, \ldots, t.$
Another new insight is that we can now see that the singular vectors of the idempotent matrix are related to the
singular vectors of the corresponding involutory matrix via a real orthogonal transformation.

%\section*{Acknowledgements}

%\bibliographystyle{plain}
%\bibliography{mybib}

\end{document}